\documentclass[12pt]{amsart}
\usepackage{amssymb}

\newtheorem{thm}{Theorem}
\newtheorem{lem}{Lemma}
\newtheorem{cor}{Corollary}
\newtheorem{fact}{Fact}
\newtheorem{prop}{Proposition}

\newcommand{\supp}{\operatorname{supp}}
\newcommand{\cl}{\operatorname{cl}}

\def\QC{\mathop{\mathrm{QC}}\nolimits}

\begin{document}

\title[Bourgain algebras, minimal
envelopes, and minimal support sets]{Bourgain algebras, minimal envelopes,\\
minimal support sets, and some applications}

\author{Carroll J.~Guillory}

\address{\hskip-\parindent Department of Mathematics\\
University of Southwestern Louisiana\\
Lafayette, LA 70504}

\curraddr{Mathematical Sciences Research Institute\\
100 Centennial Drive\\
Berkeley, CA 94720}
\email{cjg@msri.org}

\begin{abstract}
We explicitly compute certain Douglas algebras that are invariant
under both the Bourgain map and the minimal envelope map. We also
compute the Bourgain algebra and the minimal envelope of the maximal
subalgebras of a certain singly generated Douglas algebra.
\end{abstract}



\maketitle

\section{Introduction}
Let $H^\infty$ be the Banach algebra of bounded analytic functions on
the open unit disk~$D$. By considering boundary functions, we may
regard $H^\infty$ an an essentially supremum norm-closed subalgebra of
$L^\infty=L^\infty(\partial D)$. A closed subalgebra between
$H^\infty$ and $L^\infty$ is called a {\it Douglas algebra}. $H^\infty+C$ is
the smallest Douglas algebra, where $C$ is the space of continuous
functions on $\partial D$. The Chang-Marshall theorem~\cite{Ch,Ma}
says that every Douglas algebra $B$ is generated by $H^\infty$ and
complex conjugates of interpolating Blaschke products $\psi$ with
$\bar\psi\in B$. We denote by $M(B)$ the maximal ideal space
of~$B$. Then we may consider that $M(L^\infty)\subset M(B)\subset
M(H^\infty)$ and $M(L^\infty)$ is the Shilov boundary for every Douglas
algebra~$B$. For a point~$x$ in~$M(H^\infty)$, there is a representing
measure~$\mu_x$ on $M(L^\infty)$; $f(x)=\int_{M(L^\infty)}f\,d\mu_x$
for every $f\in H^\infty$. A Douglas algebra~$B$ with $A\subsetneq B$
is called a {\it minimal superalgebra} of~$A$ if there are no Douglas
algebras~$B'$ such that $A\subsetneq B'\subsetneq B$. In~\cite{GI1},
the authors proved that if~$A$ and~$B$ are Douglas algebras with
$A\subsetneq B$, then $B$ is a minimal superalgebra of~$A$ if and only
if $\supp\mu_x=\supp\mu_y$ for every $x,y\in M(A)\setminus M(B)$.

Let $Y$ be a Banach algebra with identity and let~$B$ be a closed
subalgebra of~$Y$. The {\it Bourgain algebra}~$B_b$ of~$B$ relative to~$Y$
is defined by the set of~$f$ in~$Y$ such that $\|ff_n+B\|\to0$ for
every sequence~$\{f_n\}_n$ in~$B$ with $f_n\to0$ weakly. In this note,
we consider $Y=L^\infty$ and $B$ a Douglas algebra. In~\cite{CJY},
Cima, Janson, and Yale proved that $H_b^\infty=H^\infty+C$. Gorkin,
Izuchi, and Martini~\cite{GIM} studied Bourgain algebras of Douglas
algebras and proved that $(B_b)_b=B_b$. We denote by $B_m$ the
smallest Douglas algebra which contains all minimal superalgebras
of~$B$. We call $B_m$ the {\it minimal envelope} of~$B$. We have $B\subseteq
B_b\subseteq B_m$.

A sequence $\{z_n\}_n$ in $D$ is called {\it interpolating} if for every
bounded sequence $\{a_n\}_n$ there exists $f$ in $H^\infty$ such that
$f(z_n)=a_n$ for every $n$. A Blaschke product
$$\psi(z)=\prod_{n=1}^\infty
\frac{-\bar z_n}{|z_n|}\,\frac{z-z_n}{1-\bar z_nz}\,,\quad z\in D$$
is called {\it interpolating} if its zeroes $\{z_n\}_n$ is
interpolating. For a function~$f$ in $H^\infty$ we put $Z(f)=\{x\in
M(H^\infty+C):f(x)=0\}$. For $x\in M(H^\infty)$ we denote by
$\supp\mu_x$ the support set for the representing measure $\mu_x$. For
a subset~$E$ of~$M(L^\infty)$ we denote by $\cl E$ the closure of~$E$
in~$M(L^\infty)$. Put $E=\bigcup\{\supp\mu_x: x\in M(H^\infty+C),\
|q(x)|<1\}$ where $q$ is an interpolating Blaschke product. Then
$N(\bar q)=\cl E$. A closed set $E$ in $M(L^\infty)$ is called a {\it peak
set} for a Douglas algebra $B$ if there is a $f\in B$ such that $f=1$
on~$E$ and $|f|<1$ on $M(L^\infty)\setminus E$. A closed set~$E$ in
$M(L^\infty)$ is called a {\it weak peak} for~$B$ if $E$ is the intersection
of some family of peak sets. If $E$ is a weak peak set for~$B$, then
the set
$$B_E=\{f\in L^\infty: f|_E\in B|_E\}$$
is a Douglas algebra. The sets $\supp\mu_x$ and $N(\bar q)$ are weak
peak sets for~$H^\infty$, hence $H_{\supp\mu_x}^\infty$ and $H_{N(\bar
q)}^\infty$ are Douglas algebras. A point $x\in Z(q)$, where $q$ is an
interpolating Blaschke product, is called a {\it minimal element} for
$H^\infty[\bar q]$ if there is no $y\notin M(H^\infty[\bar q])$ such
that $\supp\mu_y\subsetneq\supp\mu_x$; that is, if $y\notin
M(H^\infty[\bar q])$, then either
$\supp\mu_y\cap\supp\mu_x=\varnothing$ or
$\supp\mu_x\subseteq\supp\mu_y$. We put
$$E_x=\{\lambda\in M(H^\infty):\supp\mu_\lambda=\supp\mu_x\}$$
and call $E_x$ the {\it level set} of~$x$. For $x$,~$y$ in $M(H^\infty)$, we
put
$$\rho(x,y)=\sup\{|f(y)|: f\in H^\infty,\ \|f\|_\infty\le1,\
f(x)=0\}\,.$$
The set
$$P(x)=\{\lambda\in M(H^\infty):\rho(\lambda,x)<1\}$$
is called the {\it Gleason part} containing~$x$. We have
$P(x)\subset E_x$, and if
$x\in M(B)$ then $E_x\subset M(B)$. The map that assigns the algebra
$B_b$ to a Douglas algebra $B$ is called the {\it Bourgain map} and the map
that assigns $B_m$ to $B$ is called the {\it minimal envelope map}.

A point $x\in M(H^\infty)$ is called {\it locally thin} if there is an
interpolating Blaschke product~$q$ such that $q(x)=0$ and
$$(1-|z_{n(\alpha)}|^2)|q'(z_{n(\alpha)})|\to1$$
whenever $z_{n(\alpha)}$ is a subnet of the zero sequence $\{z_n\}_n$
of~$b$ in~$D$ converging to~$x$. We say that $q$ is {\it locally thin
at}~$x$.

We say that a minimal support point~$x$ of a Douglas algebra~$B$
corresponds to the maximal subalgebra~$A_x$ if it satisfies the
equation
$$M(A_x)=M(B)\cup P_x\quad\text{or}\quad M(A_x)=M(B)\cup E_x\,.$$
We put $\QC=(H^\infty+C)\cap\overline{(H^\infty+C)}$. Then
$$\QC=\{f\in H^\infty+C:f|_{\supp\mu_x}\text{ is constant for every }x\in
M(H^\infty+C)\}\,.$$
For $x\in M(L^\infty)$, the set
$$Q=\{y\in M(L^\infty):f(y)=f(x)\text{ for every }f\in\QC\}$$
is called a {\it $\QC$-level set}. For every $y\in M(H^\infty+C)$, there is
a $\QC$-level set $Q_y$ such that $\supp\mu_y\subset Q_y$. For any
interpolating Blaschke product $q$ set 
$$m_q=\{x\in
M(H^\infty+C):x\text{ is a minimal support point of }H^\infty[\bar
q]\}\,.$$

In this paper we give examples (Theorem~1 and Theorem~2) of Douglas
algebras that show that the Bourgain algebra of the arbitrary
intersection of Douglas algebras does not equal the intersection of
their corresponding Bourgain algebras. Theorems~1 and~2 also show that
this is true for the minimal envelopes. We also compute the Bourgain
algebra and the minimal envelope of certain maximal
subalgebras. Finally we give two applications of minimal support
points of certain interpolating Blaschke products.

All of the work in this paper was done while the author was at the
Mathematical Science Research Institute. The author thanks the
Institute for its support during this period.

\begin{lem}
Let $q$ be an interpolating Blaschke product that is of type~$G$, and
such that $Z(q)\cap P_x$ is a finite set for all $x\in Z(q)$. Then
each $x\in Z(q)$ is a locally thin point.
\end{lem}

\begin{proof}
Since $q$ is of type~$G$ and the set $Z(q)\cap P_x$ is finite there is
a factor $q_0$ of $q$ such that $Z(z_0)\cap P_x=\{x\}$. Hence
$H_{\supp\mu_x}^\infty[\bar{q_0}]$ is a minimal superalgebra of
$H_{\supp\mu_x}^\infty$. Hence
$H_{\supp\mu_x}^\infty\subset(H_{\supp\mu_x}^\infty)_b$ and so by
Theorem~5 of~\cite{MY} we have that $x$ is a locally thin point.
\end{proof}

\begin{lem}
Let $B=H_{\supp\mu_y}^\infty$, where $y$ is a trivial point. If
$B\subset B_m$, then there is an interpolating Blaschke product $q$
and an $x_0\in Z(q)$ such that $\supp\mu_{x_0}=\supp\mu_y$.
\end{lem}

\begin{proof}
If $B\subset B_m$, then by Theorem~D of~\cite{GI2} there is an
interpolating Blaschke product $q$ such that $B[\bar\psi]$ is a
minimal superalgebra of~$B$. Hence $\supp\mu_y$ is a minimal support
set of $H^\infty[\bar\psi]$. By Theorem~2 of~\cite{Gu} there is an
$x_0\in Z(\psi)$ such that $\supp\mu_{x_0}=\supp\mu_y$.
\end{proof}

It would be nice if we could prove that the converse of Lemma~2 is
true. I have been unable to do so. There are special cases when the
converse holds.

\begin{thm}
Let $q$ be a sparse interpolating Blaschke product and set
$B=H_{N(\bar q)}^\infty$. Then $B=B_b=B_m$.
\end{thm}

\begin{proof}
First we have by Theorem~C of~\cite{GI3} that there is an interpolating
Blaschke product $\psi$ (if $B\subset B_b$) such that $Z(\psi)\cap
M(B)=\{x\}$, and $M(B)=M(B[\bar\psi])\cup P_x$ for some $x$ in
$M(B)$. Since $N(\bar q)=N_0(\bar q)$, by Proposition~1 of~\cite{GI1}
we have that $B=\bigcap_{y\in Z(q)}H_{\supp\mu_y}^\infty$. We will
show that there is a $y\in Z(q)$ with $y\in P_x$. Once this is done
then we have that $P_x=E_x$ for all $x\in M(B)$ and all $\psi$ with
$Z(\psi)\cap M(B)=\{x\}$. Thus both $B_b$ and $B_m$ are generated by
the same minimal superalgebras and we have $B_b=B_m$ (see Theorem~D
of~\cite{GI3}).

To see that such a $y$ exists, note that if $x\in M(B)$ then
by~\cite{Ga}, p.~39, $\supp\mu_x\subset N(\bar q)$. By Theorem~1
of~\cite{Iz} there is a $y\in Z(q)$ such that $\supp\mu_x\subset
Q_y$. We show that $y\in P_x$. Since $q$ is sparse we have that $y$ is
unique (Lemma~4 of~\cite{Iz}), $\supp\mu_y\subset Q_y$, and
$\supp\mu_y$ is a maximal support set. Thus we have three
possibilities: (1)~$\supp\mu_x\subsetneq\supp\mu_y$,
(2)~$\supp\mu_x\cap\supp\mu_y=\varnothing$, and
(3)~$\supp\mu_x=\supp\mu_y$. We show that~(3) holds.

Suppose~(1) holds. Then by Theorem~2 of~\cite{GI2} there is an
uncountable set $\Gamma\subset Z(\psi)$ such that
$\supp\mu_m\subset\supp\mu_y$ for all $m\in\Gamma$, and
$\supp\mu_m\cap\supp\mu_n=\varnothing$ for $n,m\in\Gamma$ for $n\ne
m$. This implies that $Z(\psi)\cap M(B)$ is an infinite set, which is
a contradiction. So~(1) cannot hold.

Now suppose~(2) holds. Then $\supp\mu_x\cap\supp\mu_y=\varnothing$,
hence $\bar\psi\in H_{\supp\mu_y}^\infty$. But $Z(\psi)\cap
M(B)=\{x\}$ implies that $\bar\psi\in H_{\supp\mu_m}^\infty$ for all
$m\in Z(q)$. This implies that $\bar\psi\in B$, another
contradiction. So~(3) must hold. Since $q$ is sparse this implies that
$x\in P_y$ (or $y\in P_x$).

Since such a $y$ exists we have that $B\subset B[\bar\psi]\subset
B[\bar q]$. We're going to show that $B[\bar\psi]$ cannot be a minimal
superalgebra of~$B$. Consider the algebras $B[\bar q]$ and
$H^\infty[\bar q]$. Theorem~3.2 of~\cite{GI4} shows that any
subalgebra~$A$ of $H^\infty[\bar q]$ is the intersection of a family
of maximal subalgebras of $H^\infty[\bar q]$. Since the set
$$\bigcup_{y\in Z(q)}P_y$$
is the set of minimal support points of both $B[\bar q]$ and
$H^\infty[\bar q]$, we can use the same proof of Theorem~3.2 to show
that any subalgebra $A$ of $B[\bar q]$ is also the intersection of a
family of maximal subalgebras of $B[\bar q]$. For $m\in Z(q)$ set
$B_m=B[\bar q]\cap H_{\supp\mu_m}^\infty$. Then we have that
$$B=\bigcap_{m\in Z(q)}B_m$$
and if $H_m^\infty=H^\infty[\bar q]\cap H_{\supp\mu_m}^\infty$, then
$$H^\infty+C=\bigcap_{m\in Z(q)}H_m\,.$$
Hence if $y\in Z(q)$ such that $y\in P_x$, then
$M(B)=M(B[\bar\psi])\cup P_y$, and we get that
$$B[\bar\psi]=\bigcap_{\substack{m\in Z(q) \\
m\ne y}}B_m\,.$$

Thus if $B[\bar\psi]$ is a minimal superalgebra of~$B$, then the
algebra
$$A^*=\bigcap_{\substack{m\in Z(q)\\
m\ne y}}H_m^\infty$$
is a minimal superalgebra of $H^\infty+C$. This is impossible since
$H^\infty+C$ has no minimal subalgebra
($(H^\infty+C)_b=H^\infty+C$). Hence $B[\bar\psi]$ is not a minimal
superalgebra for~$B$. We get $B=B_b=B_m$.
\end{proof}

\begin{thm}
Let $q$ be an interpolating Blaschke product and set $\tilde
m_q=\{y\in m_q:\supp\mu_y=\supp\mu_t,\ t\text{ is a trivial
point}\}$. Set $T=\bigcap_{y\in\tilde m_q}H_{\supp\mu_y}^\infty$. Then
$T=T_b$. If $\tilde m_q=m_q$ then $T=T_m$.
\end{thm}

\begin{proof}
Suppose $T\ne T_b$. Then by Theorem~C of~\cite{GI3} there is an
interpolating Blaschke product $\psi$ such that $Z(\psi)\cap
M(T)=\{y_0\}$. By Theorem~1 of~\cite{GI1} we have
$M(T)=M(T[\bar\psi])\cup P_{y_0}$ and $y_0$ is a minimal support point
of $H^\infty[\bar\psi]$. By Proposition~6 of~\cite{MY},
$(H_{\supp\mu_{y_0}}^\infty)_b=H_{\supp\mu_{y_0}}^\infty[\bar\psi]$.
So by Theorem~5 of~\cite{MY}, $y_0$ is a locally thin point. We show
that this is not the case by showing that there is an $x_0\in\tilde
m_q$ with $\supp\mu_{x_0}=\supp\mu_{y_0}$. Suppose that
$\supp\mu_{y_0}\ne\supp\mu_{x_0}$ for all $x_0\in\tilde m_q$. Then for
each $x_0\in\tilde m_q$ one of the following can occur:
(i)~$\supp\mu_{x_0}\subsetneq\supp\mu_{y_0}$,
(ii)~$\supp\mu_{x_0}\cap\supp\mu_{y_0}=\varnothing$, or
(iii)~$\supp\mu_{y_0}\subsetneq\supp\mu_{x_0}$. We will show that none
of these can actually happen. If~(i) is true, then $|\psi(x_0)|=1$
implies that $\bar\psi\in H_{\supp\mu_{x_0}}^\infty$ since $y_0$ is a
minimal support point. If~(ii) is true for all $x_0\in\tilde m_q$ then
again $\bar\psi\in H_{\supp\mu_{x_0}}^\infty$. So if~(i) or~(ii)
happens for all $x_0\in\tilde m_q$, then $\bar\psi\in T$, which
implies that $T[\bar\psi]=T$. So $T_b=T$ here. Now if~(iii) happens,
then by Theorem~2 of~\cite{GI2} there is an uncountable set
$\Gamma\subset Z(\psi)$ such that
$\supp\mu_\alpha\subsetneq\supp\mu_{x_0}$ for all $\alpha\in\Gamma$
and $\supp\mu_\alpha\cap\supp\mu_\beta=\varnothing$ if
$\alpha\ne\beta$ and $\alpha,\beta\in\Gamma$. Since $x_0\in M(T)$,
each $\alpha\in\Gamma$ is also in $M(T)$. Thus $\Gamma\subset
Z(\psi)\cap M(T)$. This implies that $Z(\psi)\cap M(T)\ne\{y_0\}$,
which leads to a contradiction. So~(iii) cannot hold for any
$x_0\in\tilde m_q$. So if there is a $y_0\in M(T)\cap Z(\psi)$ such
that $T[\bar\psi]$ is a minimal superalgebra for~$T$, then there is an
$x_0\in\tilde m_q$ with $\supp\mu_{y_0}=\supp\mu_{x_0}$. This implies
that
$(H_{\supp\mu_{y_0}}^\infty)_b=(H_{\supp\mu_{y_0}}^\infty)_b=
H_{\supp\mu_{x_0}}^\infty[\bar\psi]$. By the remark following
Theorem~5 of~\cite{MY} we have that $x_0$ is not a locally thin
point. So no such $y_0$ exist and we have $T=T_b$.

To show that $T=T_m$ if $\tilde m_q=m_q$ we proceed as follows. By the
argument above, if $T\ne T_m$ there is a $x_0\in\tilde m_q$ such that
$\{\lambda\in M(T):|\psi(\lambda)|<1\}=E_{x_0}$, $y_0\in E_{x_0}$, and
$\{y_0\}=Z(\psi)\cap M(T)$. We show that this is a contradiction by
showing the set $Z(\psi)\cap\tilde m_q$ contains an uncountable
set. This will suffice since $\tilde m_q\subset M(T)$. Without loss of
generality we can assume that if $x,y\in\tilde m_q$, $x\ne y$, then
$\supp\mu_x\cap\supp\mu_y=\varnothing$. Let $\{z_n\}_{n=1}^\infty$ be
the zero sequence of~$q$ in~$D$. Then there is a subnet
$\{z_{n_\alpha}\}_{\alpha\in A}$ such that $z_{n_\alpha}\to x_0$. Let
$\epsilon>0$ and set
$B=A\setminus\{n_\alpha:|\psi(z_{n_\alpha})|>1-\epsilon\}$. Since
$|\psi(x_0)|<1$ and the subnet $\{|\psi(z_{n_\alpha})|\}$ converges to
$|\psi(x_0)|$, there is an $\alpha_0\in A$ such that if
$\alpha\ge\alpha_0$, we have $n_\alpha\in B$. Hence if
$\alpha_1>\alpha_0$ we have that $|\psi(z_{n_\alpha})|<1-\epsilon$ for
all $\alpha\ge\alpha_1$. Take the subnet $\{z_{n_\alpha}\}_{\alpha\in
B}$. Then on the set $\overline{\{z_{n_\alpha}\}_{\alpha\in
B}}\setminus\{z_{n_\alpha}\}_{\alpha\in B}$ we have that
$|\psi|<1-\epsilon$. Take a subsequence $\{z_{n_k}\}$ of the subnet
$\{z_{n_\alpha}\}_{\alpha\in B}$. Let $q_0$ be the factor of~$q$ with
zero sequence $\{z_{n_k}\}$. Then
$Z(q_0)=\overline{\{z_{n_k}\}}\setminus\{z_{n_k}\}$ and $|\psi(u)|<1$
for all $u\in Z(q_0)$. By Theorem~2 of~\cite{GI1} there is an
uncountable set $\Gamma\subset Z(q_0)$ such that $\Gamma\subset\tilde
m_q$. By Theorem~2 of~\cite{GI1} we can assume that if
$\alpha,\beta\in\Gamma$, $\alpha\ne\beta$, then
$\supp\mu_\alpha\cap\supp\mu_\beta=\varnothing$. For each
$\alpha\in\Gamma$ there is an $x_\alpha\in Z(\psi)$ such that
$\supp\mu_{x_\alpha}\subseteq\supp\mu_\alpha$. This implies that
$x_\alpha\in M(T)$ since $M(T)=\overline{\bigcup_{y\in\tilde
m_q}M(H_{\supp\mu_y}^\infty)}$. Thus we have that $x_\alpha\in
Z(\psi)\cap M(T)$. So the set $Z(\psi)\cap M(T)$ is uncountable. This
shows that $T=T_m$.
\end{proof}

Both $B_b$ and $B_m$ are generated by a special type of minimal
superalgebras (determined by the character of the minimal support
point).

Under certain conditions we can determine the Bourgain algebras and
the minimal envelope algebras of maximal subalgebras of a Douglas
algebra $B$ (here $B$ will have a maximal subalgebra).

\begin{thm}
Let $B$ be a Douglas algebra such that $B_b[\bar q]=B[\bar q]$ or
$B_m[\bar q]=B[\bar q]$ for some interpolating Blaschke
product~$q$. Let $A_x$ be the maximal subalgebra of $B[\bar q]$ for
which $x$ is the corresponding minimal support point of $B[\bar
q]$. Then either:
\begin{enumerate}
\item $A_x\subset(A_x)_b$ and $(A_x)_b=(A_x)_m=B[\bar q]$, or
\item $A_x=(A_x)_b$ and $(A_x)_m=B[\bar q]$.
\end{enumerate}
\end{thm}

\begin{proof}
By using Theorem~3 of~\cite{MY}, Theorem~3 of~\cite{GIM}, and
Theorems~4 and~5 of~\cite{GI3}, our hypothesis implies that $B[\bar
q]=(B[\bar q])_b=(B[\bar q])_m$.

Now let $A_x$ be any maximal subalgebra of $B[\bar q]$ associated with
the minimal support point~$x$. First we assume that $x$ is a locally
thin point (see Theorem~5 of~\cite{MY}). Then the maximal ideal space
of $B[\bar q]$ and $A_x$ are related by the equation
$$M(A_x)=M(B[\bar q])\cup P_x\,.$$
Then $B[\bar q]$ is a minimal superalgebra of $A_x$, hence $B[\bar
q]\subseteq(A_x)_b$. Now, using Theorem~3 of~\cite{MY} again, we get
\begin{align*}
B[\bar q]
&=(B[\bar q])_b\\
&=(A_x[\bar q])_b\\
&=(A_x)_b[\bar q]\\
&=(A_x)_b\quad\text{since }\bar q\in(A_x)_b\,.
\end{align*}
Similarly $B[\bar q]=(A_x)_m=(A_x)_b$ if $x$ is a locally thin
point. This proves~(i).

Now suppose that $x$ is not a locally thin point (for example, $P_x$
is not a homeomorphic disk, see~\cite{Ga}). Then $M(A_x)$ and
$M(B[\bar q])$ are related by
$$M(A_x)=M(B[\bar q])\cup E_x$$
with $P_x\subsetneq E_x$ (or $Z(q)\cap P_x$ has infinitely many
points). Using Theorem~4 of~\cite{GI3} we have
\begin{align*}
B[\bar q]
&=(B[\bar q])_m\\
&=(A_x[\bar q])_m\\
&=(A_x)_m[\bar q]\quad\text{by Theorem~4 of~\cite{Gu}}\\
&=(A_x)_m\quad\text{since }\bar q\in(A_x)_m\,.
\end{align*}

Now, for any Douglas algebra~$A$ we have that $A\subseteq A_b\subseteq
A_m$. Thus
$$A_x\subseteq(A_x)_b\subsetneq B[\bar q]=(A_x)_m\,.$$
Since $A_x$ is a maximal subalgebra of $B[\bar q]$, we have that
$A_x=(A_x)_b$ if $x$ is not locally thin. This proves~(ii).
\end{proof}

\begin{cor}
Let $q$ be any interpolating Blaschke product and set $B=H^\infty[\bar
q]$. Let $A_x$ be any maximal subalgebra of~$B$ that corresponds to
the minimal support point of~$B$. Then either
\begin{enumerate}
\item $A_x\subset(A_x)_b$ and $(A_x)_b=(A_x)_m=B$, or
\item $A_x=(A_x)_b$ and $(A_x)_m=B$.
\end{enumerate}
\end{cor}

\begin{cor}
Let $A$ be any Douglas algebra such that $A=A_b$ or $A=A_m$, and
let~$q$ be any interpolating Blaschke product such that $\bar q\notin
A$. Set $B=A[\bar q]$ and let $B_x$ be any maximal subalgebra of~$B$
corresponding to the minimal support point of~$A[\bar q]$. Then either
\begin{enumerate}
\item $B_x\subset(B_x)_b$ and $(B_x)_b=(B_x)_m=B$, or
\item $B_x=(B_x)_b$ and $(B_x)_m=B$.
\end{enumerate}
\end{cor}

\begin{thm}
Let $B$ be a Douglas algebra that has a maximal subalgebra~$A_x$,
where~$x$ is the minimal support point of~$B$ corresponding
to~$A_x$. Then $(A_x)_m=B_m$.
\end{thm}

\begin{proof}
By Theorem~4 of~\cite{GI3} we have that $(A_x)_m\subseteq B_m$ since
$A_x\subseteq B$. Since $A_x$ is a maximal subalgebra of~$B$ there is
an $x_0\in M(A)\setminus M(B)$ and an interpolating Blaschke product
$\psi_0$ such that
\begin{align*}
M(A_x)
&=M(A[\bar\psi_0])\cup E_{x_0}\\
&=M(B)\cup E_{x_0}\,.
\end{align*}
So by Theorem~D of~\cite{GI3} we have $B\subseteq(A_x)_m$. If $B_0$ is
another minimal superalgebra containing $A_x$, then there is some
$y_0\in M(A_x)$ such that $M(A_x)=M(B_0)\cup E_{y_0}$. Hence we have
that $B_0\subseteq(A_x)_m$ and $y_0\in M(B)$, otherwise
$E_{y_0}=E_{x_0}$. To show that $B_m\subseteq(A_x)_m$, let $\psi$ be
any interpolating Blaschke product such that $\bar\psi\in B_m$. Then
by Theorem~D of~\cite{GI3} we can assume that
$$\{\lambda\in M(B):|\psi(\lambda)|<1\}=E_x$$
for some $x\in M(B)$. But 
$$\{m\in M(A_x):|\psi(m)|<1\}=E_x\cup\{m\in
M(A_x):|\psi(m)|<1\}\cap E_{x_0}.
$$ The set on the right hand side is
either $E_x$ or $E_x\cup E_{x_0}$. Hence by Theorem~4 of~\cite{GI3} we
have that $\bar\psi\in(A_x)_m$. Hence $B_m\subseteq(A_x)_m$.
\end{proof}

For $(A_x)_b$ we have the following special result if we assume an
additional assumption.

\begin{thm}
Let $B$ be a Douglas algebra that has a maximal subalgebra $A_x$,
where $x$ is the minimal support point of $B$ corresponding to
$A_x$. Assume that $P_x$ is a nonhomeomorphic disk. Then
$(A_x)_b\subsetneq B_b$.
\end{thm}

\begin{proof}
Since $B\subseteq B_b$ and $(A_x)_b\subseteq B_b$, it suffices to show
that $B\nsubseteq(A_x)_b$. Since $A_x$ is a maximal subalgebra of~$B$
corresponding to~$x$, by Theorem~1 of~\cite{GI1} we have
$$M(A_x)=M(B)\cup E_x\,.$$
Note that $P_x\subset E_x$. Hence if $\psi$ is any interpolating
Blaschke product, we have by Corollary~1.5 of~\cite{GLM} the set $M(A_x)\cap
Z(\psi)\supseteq P_x\cap Z(\psi)$ is an infinite set. By Theorem~2
of~\cite{GIM}, $\bar\psi\notin(A_x)_b$. Hence $B\nsubseteq(A_x)_b$.
\end{proof}

The following two propositions on minimal support points seem to
indicate that the sets given in them are smaller in some sense than
the set in the following two well-known facts.

\begin{fact}
Let $B$ be any Douglas algebra. Then an interpolating Blaschke
product~$q$ is invertible in~$B$ if and only if $Z(q)\cap
M(B)=\varnothing$.
\end{fact}

\begin{fact}
For any Douglas algebra~$B$ we have
$$B=\bigcap_{x\in M(B)}H^\infty_{\supp\mu_x}\,.$$
\end{fact}

\begin{prop}
An interpolating Blaschke product~$q$ in invertible in a Douglas
algebra $B$ if and only if $M(B)\cap m_q=\varnothing$.
\end{prop}

\begin{proof}
By Theorem~2 of~\cite{GI1} we have that $m_q\subseteq Z(q)$, hence if
$m_q\cap M(B)\ne\varnothing$ then $\bar q\notin B$.

To prove the converse, suppose $\bar q\notin B$. Then by Fact~1,
$Z(q)\cap M(B)\ne\varnothing$. By the proof of Theorem~2
of~\cite{GI1}, there is a $y_0\in Z(q)$ such that
$\supp\mu_{y_0}\subseteq\supp\mu_x$ for any $x\in Z(q)\cap M(B)$ and
$y_0\in m_q$. Since $x\in M(B)$ we have that
$M(H^\infty_{\supp\mu_x})\subset M(B)$. Since
$M(H^\infty_{\supp\mu_x})=M(L^\infty)\cup\{\lambda\in
M(H^\infty+C):\supp\mu_\lambda\subseteq\supp\mu_x\}$ we have that
$y_0\in M(B)$. Hence $m_q\cap M(B)\ne\varnothing$.
\end{proof}

Let $B$ be any Douglas algebra and set $m_q(B)=m_q\cap M(B)$. Let
$$M_B=\bigcup\{m_q(B):q\text{ is an interpolating Blaschke product,
}\bar q\notin B\}\,.$$
Let $\Gamma(B)=\{x_\alpha\}_{\alpha\in\Lambda}$ be the family of all
minimal support points from $M_B$ such that
$\supp\mu_{x_\alpha}\cap\supp\mu_{x_\beta}=\varnothing$ if
$\alpha\ne\beta$. Then 

\begin{prop}
$\displaystyle{B=\bigcap_{x_\alpha\in\Gamma(B)}H^\infty_{\supp\mu_{x_\alpha}}}$.
\end{prop}

\begin{proof}
If $x\in M_B$, then there is an $x_\alpha\in\Gamma(B)$ such that
$\supp\mu_{x_\alpha}=\supp\mu_x$. So it suffices to show that
$$B=\bigcap_{x\in M_B}H^\infty_{\supp\mu_x}\,.$$
Set $B_0=\bigcap_{x\in M_B}H^\infty_{\supp\mu_x}$. Since $M_B\subseteq
M(B)$, by Fact~2 we have that $B\subseteq B_0$. Suppose $B\subsetneq
B_0$. Then by the Chang-Marshall Theorem~\cite{Ch,Ma} there is an
interpolating Blaschke product $q$ such that $\bar q\in B_0$ but $\bar
q\notin B$. Hence there is a $y\in M(B)$ such that $q(y)=0$. Hence
$\bar q\notin H^\infty_{\supp\mu_y}$. By Theorem~2 of~\cite{GI1} (or
Proposition~1) there is a $y_0\in Z(q)$ such that
$\supp\mu_{y_0}\subseteq\supp\mu_y$ and $y_0\in M_B$. This implies
that $\bar q\notin H^\infty_{\supp\mu_{y_0}}$, so $\bar q$ cannot be
invertible in $B_0$. Thus $B_0=B$. We are done.
\end{proof}

\end{document}